\documentclass[12pt]{article}

\usepackage{amsthm,amsmath,amsfonts}
\usepackage{latexsym}
\usepackage[T1]{fontenc}
\usepackage{a4wide}
\usepackage{times}
\usepackage{array}
\usepackage{hhline}
\numberwithin{equation}{section}

\title{Stability of 2D incompressible flows in ${\bf R}^3$}
\author{Piotr Bogus\l aw Mucha}

\begin{document}
\maketitle

\begin{center}

{Institute of Applied Mathematics and Mechanics }

{Warsaw University}

{ ul. Banacha 2, 02-097 Warszawa, Poland}

{ E-mail: p.mucha@mimuw.edu.pl}

\end{center}

\bigskip

{\bf Abstract.} We investigate the global in time stability of 
regular solutions 
with large velocity vectors to the evolutionary Navier-Stokes equation in 
${\bf R}^3$. The class of stable flows contains all two dimensional 
weak solutions. The only assumption which is required is smallness of the 
$L_2$-norm  of initial perturbation or its derivative  with 
respect to the `$z$'-coordinate in the same norm. The magnitude of 
the rest of the norm of initial datum is not restricted.

\medskip



\section{Introduction}

The paper examines global in time regular
solutions to the evolutionary Navier-Stokes equations in the whole
three dimensional space. Since the problem of regularity of weak solutions
is open and stays  one of main challenges of the present mathematics,
the study in
this area is directed mainly onto finding  special classes of solutions
with large velocity vectors or developing the theory of so-called
conditional regularity begun by Serrin \cite{ser} -- see also
\cite{KoS},\cite{sohr}.

The known theory allows to find  nontrivial classes of 
global in time solutions to the Navier-Stokes equations
with full regularity \cite{kato},\cite{sohrbook},\cite{sol}
 guaranteeing the uniqueness.
 Here we will follow this direction. One approach
is to consider  the issue of stability of known generic solutions.
The problem is well investigated for the equations in bounded domains.
Thanks to the Poincare inequality it is possible to improve
information following from the dissipation of the system
 \cite{VS},\cite{Mu1},\cite{prst}. 
 The method in these cases is just modifications
of techniques for problems for ``pure'' small data.
However in more complex cases  as in \cite{rsell},\cite{zaj} 
the idea of  proofs is not so elementary.

In the whole space the problem is more advance,
there is no Poincare inequality  and there is a need to
find extra tools. We  point two approaches to this case.
In the first one \cite{BMN},\cite{MN}, 
authors assume sufficiently largeness
of the vorticity, then using not standard theory as the Navier-Stokes 
equations they are able to obtain a large class of nontrivial
regular solutions. The second type \cite{GKMS},\cite{sawada}
 is a consequence of development of the
theory of semigroups. Thanks to it we are able to prove existence of global
in time solution for a class of initial data even with linear growth.

Our paper will show that  two dimensional  solutions are
stable in ${\bf R}^3$ under a small perturbation of the $L_2$-norm
- or even under weaker assumptions. The main idea is
based on a ``reduction''
of the original problem to the  two dimensional case.
We will follow the old idea of Olga Ladyzhenskaya \cite{lad1},\cite{lad2}.
 In main steps  the imbedding 
$H^1({\bf R}^2) \subset BMO({\bf R}^2)$
and the Marcinkiewicz-type interpolation  for general spaces will play 
an essential role. This enables a two dimensional point of view on 
the estimation for the case  in the whole ${\bf R}^3$.
The class of generic solutions can be extended, but
two dimensional  solutions with finite Dirichlet integrals (in ${\bf R}^2$) 
seem to be the best identification of this set. Note that,
because of a geometrical  structure,
the total energy of obtained solutions is infinite.
Our result shows a new large class of globally in time regular solutions
with large velocity vectors.
This way a new argument for the regularity of weak solutions
in the general case is pointed, again.

From the mechanical point of view our result can give an interpretation
 that flows with two dimensional
symmetry are stable independently from the magnitude of the constant in the 
Poincare inequality for the considered domain.

\smallskip

The subject of the paper is the evolutionary
Navier-Stokes equations in the whole three dimensional space
\begin{equation}\label{ns-org}
\begin{array}{lcr}
v_{,t}+v \cdot \nabla v- \nu \Delta v + \nabla p =F
&\mbox{ \ \ \ \ in \ } & {\bf R}^3 \times (0,T),\\
\mbox{div }v=0 &\mbox{ \ \ \ \ in \ }&  {\bf R}^3\times (0,T), \\
v|_{t=0}=v_0 &\mbox{ \ \ \ \ on \ }& {\bf R}^3,
\end{array}
\end{equation}
where $v=(v^x,v^y,v^z)$ is the sought velocity of the fluid,
$p$ its pressure, $\nu$ is the constant positive viscous coefficient,
$F$ -- represents the external data, $v_0$ is an initial datum
of the sought velocity which by $(\ref{ns-org})_2$ is required to 
satisfies the compatibility condition ${\rm div\;}v_0=0$ and 
comma `,' denotes the differentiation.

The solutions to system (\ref{ns-org}) are viewed in the form
\begin{equation}\label{v-u}
v=w+u,
\end{equation}
where $w$ is a known smooth solution and $u$ is a perturbation of it. 
Our analysis will concentrate on the system describing  function $u$.
From system (\ref{ns-org}) we obtain
\begin{equation}\label{ns-per}
\begin{array}{lcr}
u_{,t}+v \cdot \nabla u - \nu \Delta u +\nabla p = -u \cdot \nabla w
& \mbox{ \ \ \ \ in \ } &  {\bf R}^3 \times (0,T), \\
\mbox{div } u = 0 &\mbox{ \  \ \ \ in \ } &  {\bf R}^3 \times (0,T), \\
u|_{t=0}=u_0 &  \mbox{ \ \ \ \ on \ } & {\bf R}^3,
\end{array}
\end{equation}
where initial datum $u_0=v_0-w|_{t=0}$.

Let us define the class of generic solutions $w$.

\smallskip

{\it {\bf Definition.} We say that $w \in \Xi$ is a generic solution 
to system (\ref{ns-org}) iff:

$w$  -- is a smooth solution to the Navier-Stokes equations (\ref{ns-org})
with external force $F$ such that
\begin{equation}\label{cond-1}
\nabla w \in L_2((0,\infty)_t; L_2({\bf R}^2_{xy});L_\infty({\bf R}_z)).
\end{equation}
}

We distinguish one direction in ${\bf R}^3$ prescribed by the $z$-coordinate
(we denote $\bar x =(x,y,z)$).

A good identification of  class $\Xi$ is the set of two dimensional
solutions i.e.
\begin{equation}\label{i1}
w(t,x,y,z)=\tilde w(t,x,y),
\end{equation}
where $\tilde w$ is a solution to the two dimensional
Navier-Stokes equations
\begin{equation}\label{i2}
\begin{array}{lcr}
\tilde w_{,t} +\tilde w \cdot \tilde \nabla \tilde w - \nu
\tilde \Delta \tilde w + \tilde \nabla \tilde p = \tilde F
&\mbox{ \ \ in } & {\bf R}^2 \times (0,T), \\
 \tilde {\rm div}\;\tilde w =0 & \mbox{ \ \ in } &
 {\bf R}^2 \times (0,T), \\
\tilde w |_{t=0} = \tilde w_0& \mbox{ \ \ on }& {\bf R}^2
\end{array}
\end{equation}
with an analogical description as for system (\ref{ns-org}).
If $\tilde F \in L_2(0,\infty;\dot{H}^{-1}({\bf R}^2))$, then  the energy
estimate  for  solutions to (\ref{i2})  yields the inclusion
\begin{equation}\label{i3}
\tilde \nabla \tilde w \in L_2(0,\infty;L_2({\bf R}^2_{xy})).
\end{equation}
In the force-free case the description of properties of solution $\tilde w$ 
can be  better precise. The results from \cite{wieg}, \cite{AGSS} 
imply that
\begin{equation}\label{i4}
\tilde \nabla \tilde w \in L_1(0,\infty; L_\infty ({\bf R}^2)),
\end{equation}
provided suitable assumptions on the initial datum $\tilde w_0$.
The class defined by (\ref{i4}) is the kernel of the set of generic
solutions. As we will see, we will be able to ``extend'' feature  (\ref{i4})
on the whole class of functions from set $\Xi$.

The main result of the paper is the following.

\smallskip

{\it {\bf Theorem 1.} Let $w \in \Xi$.
 If $u_0 \in H^1({\bf R}^3)\cap
W^{2-2/4}_4({\bf R}^3)$. Additionally one
of two below conditions is satisfies:

\smallskip

\ (i) \ \ $||u_0||_{L_2({\bf R}^3)}$ \ \ is sufficiently small;

\smallskip

or

\smallskip

(ii) \ \ $||u_{0,z}||_{L_2({\bf R}^3)}$ \ \ is sufficiently  small,
provided
\begin{equation}\label{extra-cond}
\begin{array}{c}
||w_{,z}||_{L_5({\bf R}^3 \times (0,\infty))}
\mbox{ \ \ and \ \ }
||\nabla w_{,z}||_{L_{5/2}({\bf R}^3 \times (0,\infty))}\\[6pt]
\mbox{sufficiently small comparing to norm }
||u_0||_{H^1({\bf R}^3)\cap W^{2-1/2}_4({\bf R}^3)};
\end{array}
\end{equation}
then there exists regular unique global in time 
solution to equations (\ref{ns-org}) in  form (\ref{v-u}), where
$u$ is the solution to system (\ref{ns-per}) such that 
$u\in W^{2,1}_{4(loc)}({\bf R}^3\times (0,\infty))$ and
\begin{equation}\label{osztw}
<u>_{W^{2,1}_4({\bf R^3}\times (0,\infty))}:=
||u_t||_{L_4({\bf R}^3\times (0,\infty))}+
||\nabla^2 u||_{L_4({\bf R}^3\times (0,\infty))}\leq DATA,
\end{equation}
where $DATA$ depends on norms of initial datum $v_0$ and vector field $w$.}

\smallskip

The above result points a large class of regular global in
time solutions to the Navier-Stokes equations in ${\bf R}^3$.
From (\ref{osztw}) -- by the classical results \cite{sohrbook} --
  solutions delivered
by Theorem 1 become smooth provided smoothness of initial data.
In particular  by (\ref{i3}) we obtain that any two dimensional weak
solution -- being sufficiently smooth -- is stable in the whole
three dimensional space. Obviously smallness of a possible
perturbation depends on the magnitude of the whole norm of the perturbed flow,
however it is restricted to cases (i) or (ii), and only one of them
have to be fulfilled. Comparing to results from \cite{prst} where stability of
two dimensional flows were considered, too, our assumption (\ref{cond-1}) 
admits a larger class of generic flows. In \cite{prst} the authors required 
$\nabla w \in L_4(0,\infty;L_2(\Omega))$, additionally 
the whole $H^1({\bf R}^3)$-norm of the initial datum 
 has been assumed  to be sufficiently small.

The proof of Theorem 1 
is based on the classical energy method, however the novel
idea is to reduce the view of the nonlinear term with respect
to the geometrical structure of  given flow $w$. The energy method allows
us to obtain an information about  solutions omitting the 
influence of nonlinear convective term $v \cdot \nabla v$.  
Similar technique for  simpler versions of the presented 
 problem has been applied in  \cite{Mu1},  \cite{Mu2}.

 An alternative approach can be given by the theory of semigroups. However
 this technique requires the smallness of the whole norm and in the most
 optimal case in the three dimensions -- by  Kato's results \cite{kato}
  we ought
 to assume smallness of the $L_3$-norm of the initial datum.
  Here for any given initial norm in space $H^1({\bf R}^3)\cap
 W^{2-1/2}_4({\bf R}^3)$ we describe the required  smallness of
 the $L_2$-norm. In particular the $L_3$-norm (even any $L_{2+\epsilon}$)
 can be arbitrary large.

Theorem 1 can be stated in spaces $W^{2,1}_p$ with $p\geq 2$ defined by the norm
\begin{equation}\label{sob}
\begin{array}{c}
||u||_{W^{2,1}_p({\bf R}^3\times (0,T))}=
||u||_{L_p({\bf R}^n\times (0,T))}+<u>_{W^{2,1}_p{\bf R}^3\times (0,T))}\\[8pt]
=(\int_0^T\int_{{\bf R}^3}|u|^pdxdt)^{1/p}+
(\int_0^T\int_{{\bf R}^3}|u_t|^pdxdt)^{1/p}+
(\int_0^T \int_{{\bf R}^3}|\nabla^2_xu|^2dxdt)^{1/p}.
\end{array}
\end{equation}
The trace of a function from the $W^{2,1}_p$-space for fixed time as for 
$t=0$ belong to the Besov $W^{2-2/p}_p$-class introduced by the norm (for $p>2$)
\begin{equation}\label{sob1}
\begin{array}{c}
||u||_{W^{2-2/p}_p({\bf R}^3)}=||u||_{L_p({\bf R}^3)}+
<u>_{W^{2-2/p}_p({\bf R}^3)}
\\[8pt]
=(\int_{{\bf R}^3}|u|^pdx)^{1/p}+(\int_{{\bf R}^3}\int_{{\bf R}^3}
\frac{|\nabla_x u(x)-\nabla_x u(y)|^p}{|x-y|^{3+(2-2/p)p}}dxdy)^{1/p}.
\end{array}
\end{equation}
Such regularity would also guantantee  smoothness of 
solutions -- however in our considerations the case $p=4$ is distinguish 
and simplifies our 
calculations. Since we are interested in smooth solutions we will not 
relax regularity of initial data.

Throughout the paper we try to use the standard notation \cite{lad2}, 
\cite{sohrbook}. Generic
constants are denoted by the same letter $C$.

The paper is organized as follows. First we show a particular case of Theorem 1.
In section 3 we construct the main estimate for  case (i). Next, we show 
analogical bound for  case (ii). And in section 5 we present a proof
of global in time existence in both cases.

\section{Motivation}

The aim of this section is to show the main idea and tools of the techniques
which will be applied to  prove Theorem 1.
We analyze a special case of system (\ref{ns-org}), we consider system 
(\ref{ns-per}) for trivial solution $w\equiv 0$ with $F \equiv 0$
for  case (ii) from Theorem 1
\begin{equation}\label{m1}
\begin{array}{lcr}
v_t+v\cdot \nabla v -\nu \Delta v +\nabla p=0
& \mbox{ \ \ \ \ in } &  {\bf R}^3 \times (0,T), \\
{\rm div}\;v=0 & \mbox{ \ \ \ \ in } &  {\bf R}^3\times (0,T),\\
v|_{t=0}=v_0 & \mbox{ \ \ \ \ on } &{\bf R}^3.
\end{array}
\end{equation}
The initial datum is required to be sufficiently smooth -- in particular
 $v_0 \in H^1({\bf R}^3)$. Additionally the compatibility
 condition ${\rm div }\,v_0=0$ is assumed. 
 We want to show the following version of
 Theorem 1.

\smallskip

{\it {\bf Theorem 2.} Let $v_0\in H^1({\bf R}^3)$.
If
\begin{equation}\label{m2}
||v_{0,z}||_{L_2({\bf R}^3)} \mbox{ \ \ \ \ is sufficiently small,}
\end{equation}
then there exists  global in time regular (unique) solution to 
system (\ref{m1}).}

\smallskip

{\bf Proof.} We skip the proof of existence. Its idea is the same
as in the one presented in  section 5, where the general system 
will be considered -- see also \cite{Mu2}.
We concentrate only on a proof of the control of smallness of
  $||v_{,z}||_{L_\infty(0,\infty;L_2({\bf R}^3))}$. 
The theory guarantees us existence of weak solutions defined globally in time.
Hence, provided sufficient smoothness of them, we  find 
a suitable a priori estimate controlling smallness of mentioned quantity.

Write the energy identity for solutions to system (\ref{m1})
\begin{equation}\label{m3}
\frac{d}{dt}||v||_{L_2({\bf R}^3)}^2+2\nu||\nabla v||_{L_2({\bf R}^3)}^2=0.
\end{equation}
From (\ref{m3}) we conclude that
\begin{equation}\label{m4}
||v||_{L_\infty((0,\infty)_t;L_2({\bf R}^3))}+
||\nabla v||_{L_2( {\bf R}^3\times (0,\infty))}\leq
C||v_0||_{L_2({\bf R}^3)}.
\end{equation}

In our considerations we distinguish a one space direction, say, the 
$z$-direction. Let us differentiate  system (\ref{m1}) 
with respect to this coordinate, getting
\begin{equation}\label{m5}
\begin{array}{lcr}
v_{,zt}+v\cdot \nabla v_{,z}-\nu \Delta v_{,z}+\nabla p_{,z}=
-v_{,z}\cdot \nabla v & \mbox{ \ \ in } &  {\bf R}^3 \times (0,\infty),\\
{\rm div}\,v_{,z}=0 & \mbox{ \ \ in } &  {\bf R}^3 \times (0,\infty),\\
v_{,z}|_{t=0}=v_{0,z} & \mbox{ \ \ on } &{\bf R}^3.
\end{array}
\end{equation}
The energy method yields the following differential inequality
\begin{equation}\label{m6}
\frac{d}{dt}||v_{,z}||_{L_2({\bf R}^3)}^2+2\nu
||\nabla v_{,z}||_{L_2({\bf R}^3)}^2\leq 2\int_{{\bf R}^3}\left|
v_{,z}\cdot \nabla v v_{,z}\right| d\bar x,
\end{equation}
where $d \bar x=dxdydz$.
Hence integrating inequality (\ref{m6}) over $(0,\infty)$ we obtain
\begin{equation}\label{m7}
\begin{array}{c}
||v_{,z}||_{L_\infty((0,\infty)_t;L_2({\bf R}^3))}+
||\nabla v_{,z}||_{L_2((0,\infty)_t; L_2({\bf R}^3))}
\\[8pt]
\displaystyle
\leq C\left(  \left( \int_0^\infty dt \int_{{\bf R}^3} d\bar x
|v_{,z}|^2 |\nabla v|\right)^{1/2} +||v_{0,z}||_{L_2({\bf R}^3)}\right).
\end{array}
\end{equation}

To simplify our notation let us introduce the following quantities 
\begin{equation}\label{m8}
\begin{array}{c}
I= ||v||_{L_\infty((0,\infty)_t;L_2({\bf R}^3))}+
||\nabla v||_{L_2({\bf R}^3)},\\[8pt]
J=   ||v_{,z}||_{L_\infty((0,\infty)_t;L_2({\bf R}^3))}+
||\nabla v_{,z}||_{L_2({\bf R}^3\times (0,\infty)_t)}.
\end{array}
\end{equation}

Taking into account information from (\ref{m4}) and (\ref{m7}), assuming
finiteness of $I$ and $J$ -- we concentrate our attention only on 
finding the estimates, so above
quantities are assumed to be finite -- we conclude that
\begin{equation}\label{m9}
\nabla v \in L_2((0,\infty)_t; L_2({\bf R}^2_{xy});H^1({\bf R}_z)).
\end{equation}
From the imbedding theorem 
($H^1({\bf R}) \subset L_\infty({\bf R})$) we have
the following inequality
\begin{equation}\label{m10}
||\nabla v||_{L_2((0,\infty)_t;L_2({\bf R}^2_{xy})
;L_\infty({\bf R}_z))} \leq C I^{1/2}J^{1/2}.
\end{equation}
Employing  the interpolation inequality from the theory from 
\cite{bennet}, we get
\begin{equation}\label{m11}
\begin{array}{c}
||v_{,z}||_{L_4({\bf R}^2_{xy}\times (0,\infty)_t;L_2({\bf R}_z))}\\[8pt]
\leq C||v_{,z}||_{L_\infty(0,\infty)_t; L_2({\bf R}^2_{xy});L_2({\bf R}_z))}^{1/2}
||v_{,z}||_{L_2((0,\infty)_t;BMO({\bf R}^2_{xy});L_2({\bf R}_z))}^{1/2}\leq
CJ.
\end{array}
\end{equation}
To get the above inequality it is enough to note that
 $H^1({\bf R}^2)\subset BMO({\bf R}^2)$, then the interpolation 
relation implies
\begin{equation}\label{m12}
\begin{array}{c}
L_4({\bf R}^2_{xy}\times (0,\infty)_t;L_2({\bf R}_z))=
\\[6pt]
\left(
(L_\infty((0,\infty)_t;L_2({\bf R}^2_{xy});L_2({\bf R}_z)),
L_2((0,\infty)_t;BMO({\bf R}^2_{xy});L_2({\bf R}_z))\right)_{1/2},
\end{array}
\end{equation}
since
$\frac{1}{4}=\frac{1-1/2}{\infty}+\frac{1/2}{2}$
and  $\frac{1}{4}=\frac{1-1/2}{2}+\frac{1/2}{BMO}$ -- the constant
in (\ref{m11}) depends on  interpolation parameters.
Note that in (\ref{m11})  the classical 
Ladyzhenskaya inequality from \cite{lad1} is hidden. This inequality 
guarantees  the solvability of the regularity problem in two dimensions.

Now we are prepared to examine the first term in the r.h.s. of (\ref{m7})
which is the only difficulty in inequality (\ref{m7}). We have
\begin{equation}\label{m13}
\begin{array}{c}
\displaystyle\left[ \int_0^\infty dt \int_{{\bf R}^3}
|v_{,z}|^2 |\nabla v| d\bar x\right]^{1/2}
\\[8pt]
\displaystyle
\leq C\left[\int_0^\infty dt \int_{{\bf R}^2}
 ||\nabla v(t,x,y,\cdot)||_{L_\infty({\bf R}_z)}
 ||v_{,z}(t,x,y,\cdot)||_{L_2({\bf R}_z)}^2 dxdy \right]^{1/2}
\\[10pt]
\displaystyle
\leq C\left[||\nabla v||_{L_2({\bf R}^2_{xy}\times (0,\infty)_t;L_\infty
({\bf R}_z))}||v_{,z}||_{L_4({\bf R}^2_{xy}\times (0,\infty)_t;
L_2({\bf R}_z))}^2\right]^{1/2}
\\[8pt]
\leq C[I^{1/2}J^{1/2}J^2]^{1/2}=CI^{1/4}J^{5/4}.
\end{array}
\end{equation}
Hence using (\ref{m8}) we state inequality (\ref{m7})  as follows:
\begin{equation}\label{m14}
\begin{array}{c}
J\leq A_0 I^{1/4}J^{5/4}+J_0,
\end{array}
\end{equation}
where $J_0=C||v_{0,z}||_{L_2({\bf R}^3)}$.
If $J_0$ is so small that $A_0I^{1/4}(2J_0)^{1/4} <\frac{1}{2}$, then
from (\ref{m14}) we conclude
\begin{equation}\label{m15}
J\leq 2J_0.
\end{equation}

Thus, the smallness of initial $J_0$ implies the global in time smallness
of norms controlled by $J$ -- see (\ref{m8}).
 Here we stop the considerations for Theorem 2,
since the rest of the proof is almost the same as in the proof of the main
theorem.
Hence we claim that Theorem 2 has been proved.

\bigskip

{\bf Remark.} From the imbedding theorem in ${\bf R}^3$ we have
\begin{equation}\label{imb}
||w||_{L_6({\bf R}^3)}\leq C
||w_{,x}||_{L_2({\bf R}^3)}^{1/3}  ||w_{,y}||_{L_2({\bf R}^3)}^{1/3}
||w_{,z}||_{L_2({\bf R}^3)}^{1/3}.
\end{equation}
Smallness of $J_0$ may imply that the $L_6$-norm 
of initial datum $v_0$ is small, too.
Next, the interpolation estimate may follow the $L_3$-norm is
small, too. However it is not the case. The initial datum taken in
Theorem 2 (or in Theorem 1) may be chosen in that way the $L_3$-norm
is arbitrary large (even  $L_{2+\epsilon}$) and the only restriction is posed on
the $L_2$-norm of the derivative with respect to $z$. It has to be
  sufficiently small comparing to the magnitude of the ``rest'' of the norm of
the initial datum.

\section{Control of the $L_2$-norm}

In this part we show  the basic a priori estimate of the $L_2$-norm
of solutions to system (\ref{ns-per}). 
Precisely, we prove the estimate to part (i) of Theorem 1.

\smallskip

{\it {\bf Lemma 3.} Let $w\in \Xi$, then sufficiently smooth solutions to
(\ref{ns-per}) fulfill the following estimate
\begin{equation}\label{c1}
||u||_{L_\infty(0,\infty;L_2({\bf R}^3))}+
||\nabla u||_{L_2({\bf R}^3 \times (0,\infty)))}\leq
C||u_0||_{L_2({\bf R}^3)}.
\end{equation}
}

{\bf Proof.}
For any given $w$ fulfilling (\ref{cond-1}) and any given $\epsilon >0$ 
we are able to find a smooth function 
$Q: (0,\infty) \times {\bf R}^3 \to {\bf R}$ such that
\begin{equation}\label{c2}
||Q-|\nabla w| \,||_{L_2({\bf R}^3_{xy}\times (0,\infty)_t;L_\infty({\bf
R}_z))}\leq \epsilon
\end{equation}
and
\begin{equation}\label{c3}
Q\in L_1(0,\infty; L_\infty({\bf R}^3)).
\end{equation}

We treat system (\ref{ns-org}) as  a perturbation of a special flow $w$.
Multiplying $(\ref{ns-per})_1$ by $u$, integrating over ${\bf R}^3$, we get
\begin{equation}\label{c5}
\frac{d}{dt}||u||_{L_2({\bf R}^3)}^2+2\nu ||\nabla u||_{L_2
({\bf R}^3)}^2=-2\int_{{\bf R}^3} u \cdot \nabla w u d\bar x.
\end{equation}
Let us consider the r.h.s. of (\ref{c5}). Since the regularity or rather
vanishing conditions on function $w$ 
are too weak, we apply a trick with function $Q$ in the following way
\begin{equation*}
\left| \int_{{\bf R}^3} u \cdot \nabla w u d \bar  x \right|\leq
\int_{{\bf R}^3} |Q||u|^2 d\bar x +
\int_{{\bf R}^3}||Q-|\nabla w(t,x,y,\cdot)|\,||_{L_\infty({\bf R}_z)}
|u|^2 d \bar x.
\end{equation*}
Hence the identity (\ref{c5}) yields the following inequality
\begin{equation}\label{c7}
\begin{array}{c}
\displaystyle
\frac{d}{dt}\left[
||u||_{L_2({\bf R}^3)}^2\exp\left\{
-\int_0^t||Q||_{L_\infty({\bf R}^3)}ds\right\}\right]
+2\nu ||\nabla u||^2_{L_2({\bf R}^3)}\exp\left\{
-\int_0^t ||Q||_{L_\infty({\bf R}^3)}ds\right\}
\\[8pt]
\displaystyle
\leq C \int_{{\bf R}^3} || Q-|\nabla w|\,||_{L_\infty({\bf R}_z)}
|u|^2\exp\left\{-\int_0^t ||Q||_{L_\infty({\bf R}^3)}ds\right\}d\bar x.
\end{array}
\end{equation}

Let us introduce an auxiliary function redefining our sought function
\begin{equation}\label{c8}
{\cal U}=u \exp\left\{-\frac{1}{2} \int_0^t||Q||_{L_\infty({\bf R}^3)}
ds \right\}.
\end{equation}
And again, the same as in section 2, we introduce
\begin{equation}\label{c10}
K=||{\cal U}||_{L_\infty(0,\infty; L_2({\bf R}^3))}
+||\nabla {\cal U}||_{L_2({\bf R}^3 \times (0,\infty))}.
\end{equation}
Then inequality (\ref{c7}) can be stated as follows
\begin{equation}\label{c11}
K^2 \leq C\int_0^\infty \int_{{\bf R}^3}
||Q-|\nabla w|\,||_{L_\infty({\bf R}_z)}
 |{\cal U}|^2 d\bar x+ ||u_0||_{L_2({\bf R}^3)}^2.
\end{equation}
Assuming finiteness of $K$ the same as for (\ref{m11})-(\ref{m12}) 
we conclude
\begin{equation*}
{\cal U} \in L_4({\bf R}_{xy}^2 \times (0,\infty)_t;L_2({\bf R}_z)).
\end{equation*}
Take the first term from the r.h.s. of (\ref{c11})
\begin{equation}\label{c13}
\begin{array}{c}
\displaystyle
\int_0^\infty dt \int_{{\bf R}^3}
||Q-|\nabla w|\,||_{L_\infty({\bf R}_z)}|{\cal U}|^2 dxdydz
\leq \int_0^{\infty} dt \int_{{\bf R}^2}
||Q-|\nabla w|\,||_{L_\infty({\bf R}_z)}
||{\cal U}||^2_{L_2({\bf R}_z)}dxdy
\\[8pt]
\displaystyle
\leq C\left(\int_0^{\infty} \int_{{\bf R}^2}
||Q-|\nabla w|\,|_{L_\infty({\bf R}_z)}^2 dt dxdy\right)^{1/2}
\left(\int_0^\infty \int_{{\bf R}^2}
||{\cal U}||^4_{L_2({\bf R}_z)}dt dx dy \right)^{1/2}.
\end{array}
\end{equation}
So by (\ref{c2}) and (\ref{c13}) inequality (\ref{c11}) takes the 
following form
\begin{equation}\label{c14}
K^2 \leq C\epsilon K^2 +K_0^2
\end{equation}
with $K_0=C||u_0||_{L_2({\bf R}^3)}$.
Since  $C$ in (\ref{c14}) is an absolute constant, we can  choose 
 $\epsilon$ -- see (\ref{c2}) -- such that inequality (\ref{c14}) yields
\begin{equation}\label{c15}
K \leq 2K_0.
\end{equation}
From the definition of $K$ -- see (\ref{c10}) -- we deduce (\ref{c1}), since
by (\ref{c3}) integral $\int_0^\infty ||Q||_{L_\infty}ds$ is finite and given.
Lemma 3 is proved.

\smallskip

The obtained estimate stays independently from the magnitude of
initial datum  $K_0$. Hence
if $K_0$ is small, then $K$ is small, too. 
Lemma 3 applied to case (i) from Theorem 1 guarantees that
uniformally in time the smallness of the $L_2$-norm  is controlled.

Another advantage of Lemma 3 is that it does not require smallness
of the $L_2$-norm of initial datum $u_0$,
 hence it works in  case (ii) of Theorem 1,
too. Thus, the next section starts with information given by (\ref{c1}).

\section{Differentiation with respect to ``$z$''}

In this section we show the main estimate of the proof  of 
 the second part of Theorem 1. We prove.

\smallskip

{\it {\bf Lemma 4.} Let assumptions of Theorem 1 -- case (ii) 
with conditions (\ref{extra-cond}) be fulfilled, then sufficiently smooth
solutions to system (\ref{ns-per}) satisfy the following bound
\begin{equation}\label{d0}
||u_{,t}||_{L_\infty(0,\infty;L_2({\bf R}^3))}+
||\nabla u||_{L_2({\bf R}^3 \times (0,\infty))}\leq 
C(||u_{0,z}||_{L_2({\bf R}^3)}+\sigma ||u||_{L_2({\bf R}^3)}).
\end{equation}
where $\sigma$ describes smallness of norms mentioned in
condition (\ref{extra-cond}).}

\smallskip

{\bf Proof.} 
Differentiating system (\ref{ns-per}) with respect to the $z$-coordinate
 we get from the first (momentum) equation  the following one
\begin{equation}\label{d1}
\begin{array}{l}
u_{,zt}+v\cdot \nabla u_{,z} -\nu \Delta u_{,z} +\nabla p_{,z}= \qquad
\qquad \qquad \qquad  \qquad \qquad \qquad \qquad
\\
\qquad \qquad  \qquad \qquad
-u_{,z}\cdot \nabla u +w_{,z}\cdot \nabla u -u_{,z}\cdot \nabla w
-u \cdot \nabla w_{,z}    \mbox{ \ \ in } {\bf R}^3\times (0,\infty).
\end{array}
\end{equation}
Multiplying (\ref{d1}) by  $u_{,z}$, integrating over ${\bf R}^3$, we get
\begin{equation}\label{d2}
\begin{array}{c}
\displaystyle
\frac{d}{dt} ||u_{,z}||^2_{L_2({\bf R}^3)}
+2\nu ||\nabla u_{,z}||_{L_2({\bf R}^3)}^2\leq
C \left(
\int_{{\bf R}^3} |u_{,z} \cdot \nabla w u_{,z}| d \bar x 
+
\int_{{\bf R}^3} |u_{,z} \cdot \nabla u u_{,z}| d \bar x \right.
\\[6pt]
\displaystyle
\left.+
\int_{{\bf R}^3} |w_{,z} \cdot \nabla u u_{,z}| d \bar x +
\int_{{\bf R}^3} |u \cdot \nabla w_{,z} u_{,z}| d \bar x\right)
=I_1+I_2+I_3+I_4.
\end{array}
\end{equation}

In the case as generic solution $w$ is generated by a two dimensional 
 flow integrals $I_3$ and $I_4$ vanish and conditions (\ref{extra-cond})
 are trivially fulfilled ($\sigma$ in (\ref{d0}) is equal zero).

The same as in Lemma 3 we introduce
\begin{equation*}
{\cal U}_{,z}=u_{,z}
\exp\{ -\frac{1}{2}\int_0^t||Q||_{L_\infty({\bf R}^3)}ds\}.
\end{equation*}
Thus from (\ref{d2}) and properties of function $Q$ -- 
(\ref{c2}) and (\ref{c3}) -- we get
\begin{equation}\label{d4}
\begin{array}{c}
\displaystyle
\frac{d}{dt}||{\cal U}_{,z}||^2_{L_2({\bf R}^3)}
+2\nu ||\nabla {\cal U}_{,z}||_{L_2({\bf R}^3)}^2
\leq
C\left( \int_{{\bf R}^3} |Q-|\nabla w|| \, |{\cal U}_{,z}|^2 d\bar x
\right.
\\[8pt]
\displaystyle
\left. +
\int_{{\bf R}^3} |{\cal U}_{,z} \cdot \nabla u \,{\cal U}_{,z}| d \bar x
+\int_{{\bf R}^3} |w_{,z} \cdot \nabla {\cal U} \,{\cal U}_{,z}| d \bar x+
\int_{{\bf R}^3} |{\cal U} \cdot \nabla w_{,z} {\cal U}_{,z}| d \bar x\right)
\end{array}
\end{equation}
which leads the following inequality
\begin{equation}\label{d5}
\begin{array}{c}
\displaystyle
L \leq C\left[ \left( \int_0^\infty dt \int_{{\bf R}^3} 
|Q-|\nabla w||\,| {\cal U}_{,z}|^2 d \bar x\right)^{1/2}\right.
\\
\displaystyle
+ \left( \exp\{ \frac{1}{2}\int_0^\infty||Q||_{L_\infty({\bf R}^3)}ds\}
\int_0^\infty dt \int_{{\bf R}^3} |{\cal U}_{,z}\cdot
\nabla {\cal U} \,{\cal U}_{,z}|d \bar x \right)^{1/2}
\\
\displaystyle
\left.+
\left( \int_0^\infty dt \int_{{\bf R}^3}
 |w_{,z} \cdot \nabla {\cal U}\,
 {\cal U}_{,z}|d \bar x \right)^{1/2}+
\left( \int_0^\infty dt \int_{{\bf R}^3} {\cal U} \cdot \nabla w_{,z}
\,{\cal U}_{,z}|d \bar x \right)^{1/2} \right]
+L_0
\\
=A_1+A_2+A_3+A_4+L_0,
\end{array}
\end{equation}
where
\begin{equation}\label{d6}
L=||{\cal U}_{,z}||_{L_\infty((0,\infty)_t;L_2({\bf R}^3))}
+||\nabla {\cal U}_{,z}||_{L_2({\bf R}^3\times (0,\infty))}.
\end{equation}
and $L_0=C||u_{0,z}||_{L_2({\bf R}^3)}$.
Again, applying (\ref{c2}) and the method for estimation (\ref{c13}) 
 the first integral from the r.h.s. of (\ref{d5}) is 
bounded as follows
\begin{equation}\label{d7}
A_1 \leq C||Q-|\nabla w|\,||_{L_2({\bf R}^2_{xy}\times (0,\infty)_t;
L_\infty({\bf R}_z))}||{\cal U}_z||^2_{L_4({\bf R}^3_{xy}
 \times (0,\infty)_t; L_2({\bf R}_z))}\leq \epsilon L.
\end{equation}
To estimate $A_2$ we repeat exactly steps from section 2 -- estimation 
(\ref{m9})-(\ref{m13}) -- getting
\begin{equation}\label{d8}
\begin{array}{c}
\displaystyle
A_2 \leq \left[ \exp\{ \frac{1}{2}\int_0^\infty
 ||Q||_{L_\infty({\bf R}^3)}ds\}
 \int_0^\infty
\int_{{\bf R}^3} | {\cal U}_{,z} \cdot \nabla {\cal U}\,
{\cal U}_{,z}| d \bar x \right]^{1/2}
\\[10pt]
\displaystyle
\leq
C \exp\{\frac{1}{4} \int_0^\infty ||Q||_{L_\infty({\bf R}^3)}dx\}
K^{1/4} L^{5/4},
\end{array}
\end{equation}
where $K$ is defined by (\ref{c10}) and by our assumptions and Lemma 3 
is already given.

To estimate $A_3$ and $A_4$ we apply extra assumptions given 
by (\ref{extra-cond}) having 
\begin{equation}\label{d9}
\begin{array}{c}
A_3\leq C [||w_{,z}||_{L_5({\bf R}^3\times (0,\infty))}
||\nabla {\cal U}||_{L_2({\bf R}^3\times (0,\infty))}
||{\cal U}_{,z}||_{L_{10/3}({\bf R^3}\times (0,\infty))}]^{1/2}
\leq C \sigma K^{1/2}L^{1/2},
\\[8pt]
A_4\leq C[||\nabla w_{,z}||_{L_{5/2}({\bf R}^3\times (0,\infty))}
||{\cal U}||_{L_{10/3}({\bf R}^3\times (0,\infty))}
||{\cal U}_{,z}||_{L_{10/3}({\bf R}^3\times (0,\infty))}]^{1/2}
\leq C \sigma K^{1/2}L^{1/2},
\end{array}
\end{equation}
where we applied the parabolic imbedding into $L_{10/3}({\bf R}^3\times
(0,\infty))$.

Summing up estimates (\ref{d5})-(\ref{d9}), remembering that $Q$ is given 
and fulfills (\ref{c3}), thus the integral in the r.h.s. of (\ref{d8}) is given,
too,  we obtain the following inequality
\begin{equation}\label{d11}
L\leq \epsilon L+ 
C \exp\{\frac{1}{4} \int_0^\infty ||Q||_{L_\infty({\bf R}^3)}dx\}
K^{1/4}L^{5/4}+C\sigma K^{1/2}L^{1/2} +L_0.
\end{equation}
Smallness of $\epsilon$ -- see (\ref{c2}) -- and 
$\sigma$ -- see (\ref{extra-cond}) -- reduces (\ref{d11}) to 
the following form
\begin{equation}\label{d12}
L\leq A_1K^{1/4}L^{5/4}+\sigma K + 2L_0.
\end{equation}
Provided  $\sigma$ and $L_0$ 
such that $A_1K^{1/4}[4(L_0+\sigma K)]^{1/4}< \frac 12 $, controlling 
$K$ by Lemma 3 and (\ref{c15}), we conclude that 
\begin{equation}\label{d13}
 L \leq 4(L_0+\sigma K_0).
\end{equation}
Hence by (\ref{d13}) we get bound (\ref{d0}) guaranteeing us smallness of
the l.h.s. in this estimate. 
Lemma 4 is proved.

Now we are prepared to show estimate (\ref{osztw}) from Theorem 1.

\section{The existence}

In this section we show a proof of existence of regular 
global in time solutions
to system (\ref{ns-per}) guaranteeing by Theorem 1.
 Local in time results for these systems
follow from the standard approach and detailed proofs can be found e.g. in
\cite{Mu1},\cite{Mu2},\cite{sol}. 
Hence to obtain global in time solutions a priori estimates in a
suitable high class of regularity is required, only. 
Here it will be the $W^{2,1}_4$-space -- see (\ref{osztw}) and (\ref{sob}).
 First we  consider case (ii) which seems to be more advanced than (i).

A key element of our technique will be an application of information about global smallness of
quantity $L$ controlling by Lemma 4. A direct method seems
to be not so effective, but by the imbedding theorem  we  get a more
suitable quantity. By (\ref{imb}) we conclude
\begin{equation}\label{e1}
||u(\cdot,t)||_{L_6({\bf R}^3)}\leq C||u_{,z}(\cdot,t)||_{L_2({\bf R}^3)}^{1/3}
||u_{,x}(\cdot,t)||_{L_2({\bf R}^3)}^{1/3}
||u_{,y}(\cdot,t)||_{L_2({\bf R}^3)}^{1/3}
\end{equation}
which leads us to the following inequality
\begin{equation}\label{e2}
||u||_{L_3(0,\infty;L_6({\bf R}^3))}\leq C||u_{,z}||^{1/3}_{
L_\infty(0,\infty;L_2({\bf R}^3))}||\nabla u ||_{L_2(
{\bf R}^3\times (0,\infty))}^{2/3}.
\end{equation}
Next, let us note that the interpolation between $L_p$ spaces implies
\begin{equation}\label{e3}
\begin{array}{c}
L_4(0,\infty;L_4({\bf R}^3))=
\left( L_3(0,\infty;L_6({\bf R}^3)),
L_\infty(0,\infty;L_2({\bf R})^3)\right)_{1/4}.
\end{array}
\end{equation}
Hence remembering that the energy norm (\ref{c1}) is controlled by Lemma 3
by given data,  from (\ref{e2}) and (\ref{c1}) we obtain
\begin{equation}\label{e4}
||u||_{L_4({\bf R}^3\times (0,\infty))}\leq
C||u_{,z}||_{L_\infty(0,\infty;L_2({\bf R}^3))}^{1/3(1-1/4)},
\end{equation}
where $C$  in (\ref{e4}) contains the energy norm given by Lemma 3.
That is the reason we choose the  $W^{2,1}_4$-space to show existence of 
regular solutions to (\ref{ns-per}).
Obviously we can repeat the proof for any $W^{2,1}_p$ with general $p$
-- see \cite{Mu1}.

Now we estimate solutions in higher norms. We restate 
 problem (\ref{ns-per})   in the following  form
\begin{equation}\label{e5}
\begin{array}{lcr}
u_{,t}-\nu \Delta u +\nabla p =
-u \cdot \nabla u -w \cdot \nabla u- u \cdot \nabla w
& \mbox{in} & {\bf R}^3 \times (0,T),\\
\mbox{div }u=0 & \mbox{in} & {\bf R}^3 \times (0,T), \\
u|_{t=0}=u_0 & \mbox{on}& {\bf R}^3.
\end{array}
\end{equation}
Time $T$ -- above -- 
describes the lifespan of the maximal solution given by the local
result. Our goal is to show that we will be able to prolong this time to
$T=\infty$ at the end of our analysis.

By the classical results \cite{giso},\cite{Mu2},\cite{sol} for the
Stokes system in the whole space (the l.h.s. of (\ref{e5}))
  the following $L_p$-Schauder type estimate for solutions to (\ref{e5})
is known
\begin{equation}\label{e6}
\begin{array}{c}
||u_{,t}||_{L_p({\bf R}^3\times (0,T))}
+||\nabla^2u||_{L_p({\bf R}^3\times (0,T))}
\\[8pt]
\leq
C\left(||{\rm r.h.s. of} (\ref{e5})_1||_{L_p({\bf R}^3\times (0,T))}+
||u_0||_{W^{2-2/p}_p({\bf R}^3)}\right),
\end{array}
\end{equation}
where $C$ does not depend on $T$, so we can put $T=\infty$ in
estimate (\ref{e6}). In our case we consider bound (\ref{e6}) for $p=4$.

To apply estimate (\ref{e6}) there is a need  to find bound on 
the r.h.s. of (\ref{e5}) in the $L_4$-norm.

The imbedding theorem \cite[Chap. 11]{bin}  yields the following inclusions
\begin{equation*}
W^{2,1}_4 ({\bf R}^3 \times (0,T))  \subset L_{12}({\bf R}^3\times
(0,T)),\quad 
\nabla  W^{2,1}_4 ({\bf R}^3 \times (0,T))  \subset L_{6}({\bf R}^3\times
(0,T)),
\end{equation*}
moreover
there exists a function $c(\cdot)$ such that $c(\sigma) \to \infty$ as $\sigma
\to 0$ and
\begin{equation}\label{e9}
\begin{array}{c}
||u||_{L_{12}({\bf R}^3\times (0,T))}+
||\nabla u ||_{L_6({\bf R}^3 \times (0,T))}
\leq \sigma <u>_{W^{2,1}_4
({\bf R}^3 \times (0,T))}+c(\sigma)||u||_{L_4({\bf R}^3 \times (0,T))},
\end{array}
\end{equation}
where $<\cdot>_{W^{2,1}_4}$ denotes the main seminorm of space 
$W^{2,1}_4({\bf R}^3\times (0,T))$ --  see (\ref{sob}).

Applying estimate (\ref{e9}) to terms of the r.h.s. of (\ref{e6})
we get
\begin{equation}\label{e10}
\begin{array}{c}
||u\cdot \nabla u||_{L_4({\bf R}^3\times (0,T))}\leq
C||u||_{L_{12}({\bf R}^3\times (0,T))}
||\nabla u||_{L_6({\bf R}^3\times (0,T))}
\\[7pt]
\leq
\sigma^2<u>^2_{W^{2,1}_4({\bf R}^3\times (0,T))}+
c(\sigma)||u||^2_{L_4({\bf R}^3\times (0,T))}
\end{array}
\end{equation}
and
\begin{equation}\label{e11}
\begin{array}{c}
||w\cdot \nabla u||_{L_4({\bf R}^3\times (0,T))}\leq
\sigma <u>_{W^{2,1}_4({\bf R}^3\times (0,T))}+
c(\sigma,||w||_{L_\infty({\bf R}^3 \times (0,T))})
||u||_{L_4({\bf R}^3\times (0,T))},
\\[8pt]
||u\cdot \nabla w||_{L_4({\bf R}^3\times (0,T))}\leq
C||\nabla w||_{L_\infty({\bf R}^3\times (0,T))}
||u||_{L_4({\bf R}^3\times (0,T))}.
\end{array}
\end{equation}
Inserting (\ref{e10}) and (\ref{e11})
 to estimate (\ref{e6}), remembering about (\ref{e4}), we obtain
\begin{equation}\label{e12}
\begin{array}{c}
<u>_{W^{2,1}_4({\bf R}^3\times (0,T))}\leq
\sigma <u>_{W^{2,1}_4({\bf R}^3\times (0,T))}+
\sigma^2 <u>^2_{W^{2,1}_4({\bf R}^3\times (0,T))}
\\[8pt]
+c(\sigma,||w||_{W^1_\infty({\bf R}^3\times (0,\infty))})
||u||_{L_4({\bf R}^3\times (0,\infty))}
+C<u_{0}>_{W^{2-1/2}_4({\bf R}^3)}.
\end{array}
\end{equation}

Provided smallness of $\sigma$, remembering that the $L_4$-norm of $u$ by 
bound (\ref{e4}) 
is sufficiently small by (\ref{e12}), from (\ref{e12}) we obtain
\begin{equation}\label{e13}
<u>_{W^{2,1}_4({\bf R}^3\times (0,T))}\leq DATA.
\end{equation}
Note that to obtain (\ref{e13}) 
 smallness of $<u_0>_{W^{2-1/2}_4({\bf R}^3)}$ 
is not required, the only condition on
\begin{equation}\label{e14}
(1-\sigma)^2>4\sigma^2\left[   c(\sigma,
||w||_{W^1_\infty({\bf R}^3\times (0,\infty))})
||u||_{L_4({\bf R}^3\times (0,\infty))}
+C<u_{0}>_{W^{2-1/2}_4({\bf R}^3)}\right].
\end{equation}
But the choice of $\sigma$ is arbitral, additionally it  prescribes
the smallness of the $L_4$-norm of $u$ by (\ref{e4}), thus the r.h.s. 
of (\ref{e14}) can be arbitrary small. 

$DATA$ in (\ref{e13}) are bounded by  all given data,  in 
general case it may not be small.

However, first of all the r.h.s. of (\ref{e13}) does not depend on $T$, hence 
we are able to  extend our estimate on $T=\infty$, getting the desired global
in time solutions with sufficiently high regularity guaranteeing the smoothness.
Thus, we proved  case (ii) for Theorem 1.

\smallskip

Let us briefly look on  case (i). This part of Theorem 1 is similar to  case (ii),
so we point a reduction of this case to the first considered one. 

From Lemma 3 and the parabolic imbedding we immediately obtain smallness of the
$L_{10/3}$-norm, i.e.
\begin{equation}\label{z1}
||u||_{L_{10/3}({\bf R}^3 \times (0,\infty))}\leq C||u_0||_{L_2({\bf R}^3)}.
\end{equation}
Additionally the theory from \cite[Chap. 18]{bin}
 guarantees us an analogical estimate 
(\ref{e9}), but with the $L_{10/3}$-norm, i.e. there exists a function 
$c(\sigma)\to \infty$ and $\sigma \to 0$ such that
$$
||u||_{L_{12}({\bf R}^3 \times (0,\infty))}+
||\nabla u ||_{L_6({\bf R}^3 \times (0,\infty))}\leq
\sigma  <u>_{W^{2,1}_4({\bf R}^3 \times (0,\infty))}+
c(\sigma)||u||_{L_{10/3}({\bf R}^3 \times (0,\infty))}.
$$
Thus, remembering (\ref{z1}),
the whole estimation (\ref{e6})-(\ref{e13}) is almost the same. 
Concluding in a similar way we are able to show 
\begin{equation}\label{z2}
<u>_{W^{2,1}_4({\bf R}^3 \times (0,T))}\leq DATA.
\end{equation}
The same for (\ref{e13}) we can obtain bound (\ref{z2}) on time 
interval $(0,\infty)$.

The proof of  Theorem 1 is done.

{\footnotesize {\bf Acknowledgments.} The author would like to thank 
Konstantin Pileckas and
Wojciech Zaj\c aczkowski for helpful discussions.
The work has been supported by
Polish KBN grant No. 1 P03A 021 30 and by ECFP6 M.Curie ToK program SPADE2,
MTKD-CT-2004-014508 and SPB-M.}


\medskip


\bigskip

\begin{thebibliography}{99}

\bibitem{AGSS}
Ch. Amrouche; V. Girault; M.E. Schonbek; T.P. Schonbek, 
 Pointwise decay of solutions and of higher derivatives to 
Navier-Stokes equations.  SIAM J. Math. Anal.  31  (2000), 
 no. 4, 740--753

\bibitem{BMN}
A. Babin; A. Mahalov; B. Nicolaenko,  Global regularity 
of 3D rotating Navier-Stokes equations for resonant domains. 
 Indiana Univ. Math. J.  48  (1999),  no. 3, 1133--1176.

\bibitem{VS} H. Beirão da Veiga; P. Secchi,  $L\sp p$-stability 
for the strong solutions of the Navier-Stokes equations in 
the whole space.  Arch. Rational Mech. Anal.  98  (1987),  no. 1, 65--69.

\bibitem{bennet} C. Bennett; R. Sharpley, Interpolation of operators. 
Pure and Applied Mathematics, 129. Academic Press, Inc., Boston, MA, 1988.

\bibitem{bin} O. V. Besov, V. P. Ilin, and S. M. Nikolskij, 
Integral Function Representation and Imbedding Theorem, Moskow, 1975.

\bibitem{giso} Y. Giga; H. Sohr, Abstract $L\sp p$ estimates for 
the Cauchy problem with applications to the Navier-Stokes equations in 
exterior domains.  J. Funct. Anal.  102  (1991),  no. 1, 72--94.

\bibitem{GKMS} Y. Giga; K. Inui; A. Mahalov; S. Matsui, 
Navier-Stokes equations in a rotating frame in ${\bf R}\sp 3$ with 
initial data nondecreasing at infinity.  Hokkaido Math. J.  
35  (2006),  no. 2, 321--364.

\bibitem{sawada} M. Hieber; O. Sawada,  The Navier-Stokes equations 
in ${\bf R}\sp n$ with linearly growing initial data.  
Arch. Ration. Mech. Anal.  175  (2005),  no. 2, 269--285.

\bibitem{kato} T. Kato, Strong $L^p$-solutions of the Navier-Stokes equations 
in $R^m$, with application to weak solutions,  Math. Z. 187 (1984),
 471--480.

\bibitem{KoS} H. Kozono; H. Sohr,  Remark on uniqueness of
 weak solutions to the Navier-Stokes equations.  Analysis  16  
 (1996),  no. 3, 255--271.

\bibitem{lad1} O. A. Ladyzhenskaya, Solution `in the large' of 
the non-stationary 
boundary value problem for the Navier---Stokes system with two 
space variables,  Comm. Pure Appl. Math. 12 (1959), 427--433.

\bibitem{lad2} O. A. Ladyzhenskaya, The Mathematical Theory of Viscous 
Incompressible Flow, Gordon and Breach, New York, 1969.

\bibitem{MN} A.S. Makhalov; V.P. Nikolaenko,  Global solvability of 
three-dimensional Navier-Stokes equations with uniformly high 
initial vorticity. (Russian)  Uspekhi Mat. Nauk  58  (2003),  
no. 2(350), 79--110;  translation in  Russian Math. Surveys 
 58  (2003),  no. 2, 287--318. 

\bibitem{Mu1} P.B. Mucha,  Stability of nontrivial solutions 
of the Navier-Stokes system on the three dimensional torus.  
J. Differential Equations  172  (2001),  no. 2, 359--375.

\bibitem{Mu2} P.B. Mucha,  Stability of constant solutions 
to the Navier-Stokes system in ${\bf R}\sp 3$.  Appl. Math. 
(Warsaw)  28  (2001),  no. 3, 301--310. 

\bibitem{prst} G. Ponce, R. Racke, T. C. Sideris, and E. S. Titi, 
Global stability 
of large solutions to the $3{\rm D}$ Navier-Stokes equations,  
Comm. Math. Phys. 159 (1994), 329--341.

\bibitem{rsell} G. Raugel, G.R. Sell, Navier--Stokes equations on thin 
$3{\rm D}$ domains I: global attractors and global regularity of
 solutions, J. Amer. Math. Soc. 6 (1993) 503--568.

\bibitem{ser} J. Serrin, The initial value problem for the 
Navier-Stokes equations.  1963  Nonlinear Problems 
(Proc. Sympos., Madison, Wis.  pp. 69--98.

\bibitem{sohr} H. Sohr,  A regularity class for the Navier-Stokes equations 
in Lorentz spaces. Dedicated to the memory of Tosio Kato.  
J. Evol. Equ.  1  (2001),  no. 4, 441--467.

\bibitem{sohrbook} H. Sohr,  The Navier-Stokes equations. 
An elementary functional analytic approach. 
 Birkhauser Verlag, Basel, 2001.

\bibitem{sol} V.A. Solonnikov,  Estimates of the solutions of the 
nonstationary Navier-Stokes system. (Russian) 
Boundary value problems of mathematical physics and related 
questions in the theory of functions, 7.  Zap. Nau\v cn. 
Sem. Leningrad. Otdel. Mat. Inst. Steklov. (LOMI)  38  (1973), 153--231.

\bibitem{wieg}  M. Wiegner,  Decay results for weak solutions 
of the Navier-Stokes 
equations on $R\sp n$.  J. London Math. Soc. (2)  35  (1987), 
 no. 2, 303--313.

\bibitem{zaj} W. M. Zajaczkowski,  Global special regular solutions 
to the Navier--Stokes equations in a cylindrical domain under 
boundary slip conditions,  Gakuto Series in Math. 21 (2004),

\end{thebibliography}
\end{document}